\documentclass{article}

\title{On the Distribution of Witnesses in the Miller-Rabin Test}

\author{Matt Kownacki}
\usepackage{amsbsy,amssymb,amscd,amsfonts,latexsym,amstext,delarray, amsmath,graphicx,color,caption}
\usepackage{graphicx}
\usepackage{amsmath, amsthm, amssymb, amsfonts}
\usepackage{graphicx}
\usepackage{listings}
\usepackage{indentfirst}
\usepackage{cite}
\usepackage{euscript}
\usepackage{multirow}
\usepackage{etex, pictexwd,dcpic}

\newtheorem {theorem}    {Theorem}
\newtheorem{defn}{Definition}

\newtheorem{lemma}[theorem]{Lemma}

\theoremstyle{remark}

\def\Z{{\mathbb Z}}

\newcommand{\T}{\ \ \ \ \ \ \ }
	
\begin{document}
\maketitle
\textbf{Abstract:}	 We  show that the set of normalized Miller-Rabin witnesses becomes equidistributed in the unit interval. This will be done by exhibiting cancellation in certain exponential sums. 

\tableofcontents
\pagebreak

	\section{Introduction and Notation}	
	For convenience, the following notation will be put to use. For a set $S$, $\#S$ will denote the number of elements of $S$.The greatest common divisor of two integers $a$ and $b$ will be represented as $(a,b)$.We will denote the group of units modulo $n$ as $(\mathbb{Z}/n\mathbb{Z})^*$.   A function $f(n)$ is said to be $o(g(n))$ if $\lim\limits_{n \to \infty} \frac{f(n)}{g(n)} = 0$. Likewise, a function $f(n)$ is said to be $O(g(n))$ if $|f(n)| \leq c|g(n)|$ for some constant $c$. The function $e(x)$ is an exponential function to be defined as $e(x)= e^{2 \pi i x}$.

	\begin{defn}
		Let $n$ be an odd integer and write $n-1 = d 2^{s} $ with $d$ odd. Then an integer $a$, $1<a<n$,  is a Miller-Rabin witness (for the compositeness of $n$) if the following conditions hold   \\
		(1) $(a,n)=1$  \\
		(2) $a^{d} \not \equiv 1 \mod{n}$ \\
		(3) For all integers $j$ with $ 0 \leq j<s, \  a^{d2^{i}} \not \equiv -1 \mod{n}  $ \\ \end{defn}
	
	Given $n$, let $W(n)$ denote the set of such witnesses. There are two theorems of note proved by Miller and Rabin respectively.
	
	\textbf{Theorem [Miller] \cite{MR0480295}} \\ Let $n$ be odd and composite. Assuming the Generalized Riemann Hypothesis, then $$\min W(n) = O(\log(n)^2)$$
	
	The specific constant was later proved to be 2 by Erich Bach \cite{MR1023756}, so the least witness would be no larger than $2 \log{n}^2 $, assuming GRH. 
	
	\textbf{Theorem [Rabin] \cite{MR566880}}  \\ Let $n$ be odd and composite, then $\#W(n)$ obeys the following bound $$ \#W(n)>\frac{3(n-1)}{4}$$
	
	This theorem allowed Rabin to alter the deterministic version into a probabilistic version of the test, the Miller-Rabin primality test. The purpose of this paper is to study the (normalized) distribution of $W(n)$ in $\Z/n\Z$. The main result being the following. 		 
	\begin{theorem}
		As $n \to \infty$ along odd composite numbers, then the normalized witness set becomes equidistributed in the unit interval. 
	\end{theorem}
	\begin{flushleft}
		By this we mean : $\forall [a, b] \subset [0,1]$  $$ \frac{\#\{\frac{W(n)}{n}  \cap [a,b]  \}}{\#W(n)} \to b-a \ \ , as \ n \to \infty$$
		To illustrate this, see Figure 1. The proof is elementary, the main ingredient being reduction of certain exponential sums into Gauss sums. 
		\end{flushleft}
	\begin{figure}[h]
		\includegraphics[width=1\textwidth]%
		{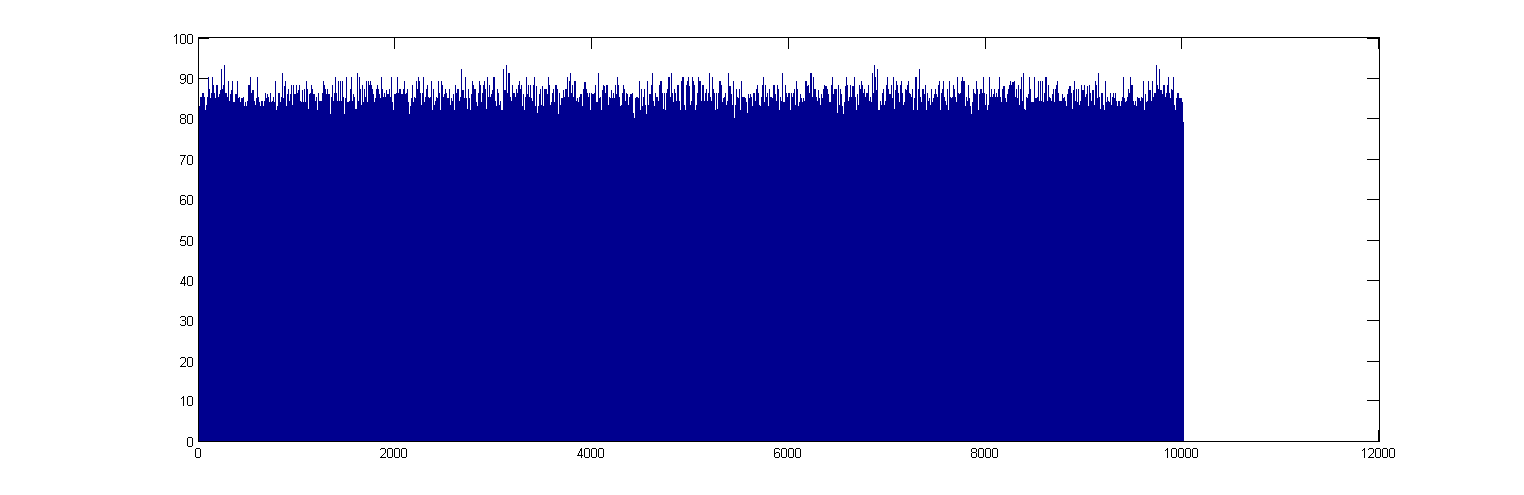}
		\\
		\caption{Witnesses of n=1056331 per interval in intervals of length $\frac{n}{10000}$.}
	\end{figure}

	\section{Proof of Theorem 1}
	\begin{flushleft}
		
	\textit{Proof.}	By Rabin's bound for $\#W(n)$ and Weyl's Criterion, it suffices to show that for fixed $k \not =0, k \in \mathbb{Z}$ that \begin{enumerate}
			\item[($\star$)] $ \T \T \T \T \T \T S= \sum_{w \in W(n)} e(\frac{kw}{n}) =o(n)$
		\end{enumerate}
		
	\end{flushleft}
	\pagebreak

	\begin{flushleft}	
		 Let $\overline{W(n)} = \{0,1,2,...,n-1\} \setminus W(n)$ be the set of non-witnesses, and define
		 $$\overline{S} = \sum_{w \in \overline{W(n)}} e(\frac{kw}{n})  $$
		 so that $S+ \overline{S} =0$. Then since $n \not | \  k$ , for $k \not = 0$, ($\star$) is equivalent to 
		\begin{enumerate}
			\item[($\dagger$)] $ \T \T \T \T \T \T | \overline{S} |  =o(n)$
		\end{enumerate}
		\begin{flushleft}
			$\overline{W(n)}$ can be partitioned based on its  membership conditions as follows: 
			
			\begin{center}
				\begin{flushleft}
					$\overline{W_{1}(n)}= \{w \in \overline{W(n)}| (w,n) >1\}$. \\
					
					$\overline{W_{2}(n)}= \{w \in \overline{W(n)}| w^d \equiv 1\mod{n} \ \}$. \\
					
					$\overline{W_{3}(n)}= \{w \in \overline{W(n)}|\exists j<s, w^{2^j d} \equiv -1 \mod{n} \}$. \\
				\end{flushleft}
			\end{center}
			
		\end{flushleft}
		Hence, $\overline{W(n)} =\overline{W_{1}(n)}\bigsqcup\overline{W_{2}(n)}\bigsqcup \overline{W_{3}(n)} $ and  ($\dagger$) follows from showing
		$$S_j= \sum_{w \in \overline{W_j(n)}} e(\frac{kw}{n}) = o(n) \ \ \ \ \ \   j=1,2,3$$\\
	
	\end{flushleft}

	\subsection{Estimation of $S_1$}
	\begin{lemma}  For each fixed $k\not = 0$, $|S_1| =O_k(1)$ as $n \to \infty$. 
	\begin{proof}
		For $S_{1}$  we can represent the sum over those $w$ as $$ \sum_{(w,n)>1} e(k\frac{w}{n})=  -\sum_{(w,n)=1} e(k\frac{w}{n})$$ \T  Then upon  a M\"{o}bius inversion  \cite{MR1395371} of this sum we arrive at:
		$$\sum_{(w,n)=1} e(k\frac{w}{n})= \sum_{s|(n,k)}s \mu(\frac{n}{s})$$ \T As this sum is just a divisor sum,and that $(n,k) \leq k$, we find that for any $\epsilon >0 $, $$ \left| S_1 \right| = \left |\sum_{s|(n,k)}s \mu(\frac{n}{s}) \right| \leq |\sum_{s|(n,k)}s|< k^{1+\epsilon} $$
		\T	Thus  $S_{1}$ is of order $O_k(1) $  as $n \to \infty$. \\ \ \\
	\end{proof}
	\end{lemma}
	
	\subsection{Cancellation Lemma}
\begin{flushleft}		
The following lemma will be applied to the estimation of $S_2$ and $S_3$.
		 \end{flushleft}
		\begin{lemma}
	
		\textit{Let $ \alpha,\  n \in \mathbb{N}$ and let $b$ be an  element in $(\mathbb{Z} / n\mathbb{Z})^*$. Let $W'= \{w \in (\mathbb{Z} / n\mathbb{Z})^*$ $|$  $w^{\alpha} \equiv b \mod{n} \}$. Fix  $k \not = 0 \in \mathbb{Z}$  , and define the sum} $$ S'= \sum_{w \in W'} e(k\frac{w}{n})$$  \textit{then} $ |S'| = O_k(\sqrt{n}) $ \textit{as} $n \to \infty$. 
	
	\begin{proof}
		Let $\overline{b}$ be the inverse of $b  \in (\mathbb{Z} / n\mathbb{Z})^*$ and consider $ \frac{1}{\phi(n)} \sum_{\chi}\chi(\overline{b} w ^ \alpha)$ , a sum over Dirichlet characters modulo $n$ , note that	$$\frac{1}{\phi(n)} \sum_{\chi}\chi(\overline{b} w ^ \alpha)=  \left\{
		\begin{array}{lr}
			1 & \T w \in W'\\
			0 & \ \T  otherwise
		\end{array}
		\right. $$
		
		 Insert this into the sum and interchange the order of summation to obtain
		
		$$ S'= \left| \sum_{w ^\alpha \equiv b \mod{n} } e(k\frac{w}{n}) \right| =\left | \sum_{w ^\alpha \equiv b \mod{n} } e(k\frac{w}{n}) \frac{1}{\phi(n)} \sum_{\chi}\chi(\overline{b} w ^ \alpha ) \right |=$$ $$\left | \frac{1}{\phi(n)} \sum_{\chi} \chi(\overline{b}) \sum_{w \mod{n}} e(k\frac{w}{n}) \chi^\alpha( w  ) \right | \leq \frac{1}{\phi(n)}\sum_{\chi} \left| \sum_{w \mod n} \chi^\alpha(w) e(k\frac{w}{n}) \right | $$
		\T	
	For each $\chi$, $\chi^\alpha$ could be the trivial character mod n or not. If it is the trivial character then the inside sum breaks down into a Ramanujan sum, and is estimated as in Lemma 2.  If it is nontrivial, note that $\sum_{w \mod n} \chi^\alpha(w) e(k\frac{w}{n}) $ is a type of Gauss sum, and it is a known fact \cite{MR2312337} that for primitive characters $$ \left |\sum_{w \mod n} \chi^\alpha(w) e(\frac{kw}{n}) \right|  \leq \sqrt{n} $$  
	
	If $\chi^\alpha$ is imprimitive with conductor $q$ , writing $n=ql$,  there are two cases to be handled. If $l \ \not |  \ k$, then the sum is zero. If $l \ | \ k$  then  $$\left | \sum_{w \mod{n}}\chi^\alpha(w)e(k\frac{w}{n}) \right | \leq
	l \sqrt{q}  $$  
	We have that $\chi^\alpha$ is induced by a character $\chi_1$ which is primitive modulo $q$. Upon writing $w=qj + r$ with $r \mod{q}$ and $j \mod{l}$ we have that  $$\sum_{w \mod{n}}\chi^\alpha(w)e(k\frac{w}{n}) =  \sum_{j\mod{l} } \ \  \sum_{r\mod{q}}\chi_1(qj+r)e(k\frac{qj + r}{ql})  $$
	
	Denote the inside sum by $$S''=\sum_{r\mod{q}}\chi_1(qj+r)e(k\frac{qj + r}{ql})$$ 
	
	To handle the first case, that $l \not | \  k$ , multiply $S''$ by $e(\frac{k}{l})$ and note that $e(\frac{k}{l})S''=S''=0$.
	For the second case, we have $k=k'l$ , and 
	$$ S''=\sum_{r\mod{q}}\chi_1(qj+r)e(k\frac{qj + r}{ql}) = \sum_{r \mod{q}}\chi_1(r)e(\frac{k}{l}\frac{r}{q} ) =\sum_{r \mod{q}}\chi_1(r)e(k'\frac{r}{q} ) $$ 
	
	Now $S''$ can be handled as above(as $\chi_1$ is primitive modulo $q$ ), and it is summed exactly $l$ times. As we have $l \ | \ k$ and $q \ | \ n$, we have $ l \leq k$ and $q \leq n$, so $\sum_{w \mod n} \chi^\alpha(w) e(k\frac{w}{n}) =O_k(\sqrt{n})$.  
	
	It follows that $|S'| = O_k(\sqrt{n}) $ as $n \to \infty$.

	\end{proof}
	\end{lemma}
	\subsection{Estimation of $S_2$ and $S_3$}
\begin{lemma} For each fixed $k\not = 0$, $|S_2| = O_k( \sqrt{n}) $ as $n \to \infty$.
	\begin{proof}
		$S_2$ is a sum over the set $\overline{W_2(n)}= \{w \in \overline{W(n)}| w^d \equiv 1\mod{n}\}$. As $\overline{W_2(n)}$ is of the type $W'$, then we can apply Lemma 3. 
	\end{proof}
	
\end{lemma}
 \begin{lemma} For each fixed $k\not = 0$, $|S_3| = O_k(\sqrt{n} \log{n}) $ as $n \to \infty$.
	
	\begin{proof}
		$S_3$ is the sum over the set 	$\overline{W_{3}(n)}= \{w \in \overline{W(n)}|\exists j<s, w^{2^j d} \equiv -1 \mod{n}\}$. We can write $S_3$ as: 
		$$S_{3}= \sum_{j=0}^{s-1} \ \ \sum_{w^{2^j d} \equiv -1 \mod{n} } e(k\frac{w}{n})$$
		Applying the cancellation lemma to the innermost sum and then bounding the outside sum by $\log{n}$ yields that $|S_3| = O_k(\sqrt{n} \log{n})$ as $n \to \infty$. 
		
	\end{proof}
	\end{lemma}
	
	\subsection{Proof of Theorem 1}

		As previously stated, the statement ($\star$)  is equivalent to $$ \left |S \right |=  \left |\sum_{w \in \overline{W(n)}} e(\frac{kw}{n}) \right | \leq \left |S_1 \right | +  \left | S_2 \right | + \left | S_3 \right | = o(n) $$ (with each sum shown individually to be $o(n)$ in Lemmas 2,4, and 5 respectively). The main result then follows from Weyl's criterion.

	\section{Acknowledgments}
	\T	The author would like to express an enormous amount of gratitude to his advisor,  Alex Kontorovich, for his wealth of experience and indispensable wisdom, Katie McKeon, whose advice and assistance made this possible, and numerous other students, staff and faculty at the DIMACS REU for their encouragement and support.

	\bibliographystyle{plain}
	\bibliography{witness}

\begin{thebibliography}{1}

\bibitem{MR1023756}
Eric Bach.
\newblock Explicit bounds for primality testing and related problems.
\newblock {\em J. Math. Comp.}, 55(191):355--380, 1990.

\bibitem{MR2312337}
Cohen. Henri.
\newblock {\em Number Theory Vol. I. Tools and Diophantine Equations}, volume
  239 of {\em Graduate Texts in Mathematics}.
\newblock Springer-Verlag, New York, 2007.

\bibitem{MR0480295}
Gary~L. Miller.
\newblock Riemann's hypothesis and tests for primality.
\newblock {\em J. Comput. System Sci.}, 13(3):300--317, 1976.
\newblock Working papers presented at the ACM-SIGACT Symposium on the Theory of
  Computing (Albuquerque, N.M., 1975).

\bibitem{MR1395371}
Melvyn~B. Nathanson.
\newblock {\em Additive number theory}, volume 164 of {\em Graduate Texts in
  Mathematics}.
\newblock Springer-Verlag, New York, 1996.
\newblock The classical bases.

\bibitem{MR566880}
Michael~O. Rabin.
\newblock Probabilistic algorithm for testing primality.
\newblock {\em J. Number Theory}, 12(1):128--138, 1980.

\end{thebibliography}

	\pagebreak
	\appendix{}
\end{document}